\documentclass[12pt]{article}
\usepackage{graphics}
\usepackage{graphicx}
\usepackage{amsmath,amsxtra,amssymb,latexsym, amscd}
\usepackage[left=3cm,right=2.5cm,top=2.5cm,bottom=2.5cm]{geometry}

\newtheorem{thm}{Theorem}

\newtheorem{cor}[thm]{Corollary}
\newtheorem{lem}[thm]{Lemma}

\begin{document}
\newcommand{\lon}{\longrightarrow}
\newcommand{\lom}{\longmapsto}
\newcommand{\dis}{\displaystyle}

\title{\bf On the Four Vertex Theorem in planes with radial density $e^{\varphi(r)}$}

\author{\bf Doan The Hieu and Tran Le Nam \\
\sl  College of Education, Hue University\\
 34 Le Loi, Hue, Vietnam\\
\sl deltic@dng.vnn.vn; dthehieu@yahoo.com}
\maketitle

\begin{abstract}
 It is showed that  in a plane with a radial density the Four Vertex Theorem holds  for the class of all  simple closed curves if and only if the density is constant. But for the class of simple closed curves that are invariant under a rotation about the origin, the Four Vertex Theorem holds for every radial density.
\end{abstract}
\noindent {\bf AMS Subject Classification (2000):}  {Primary 53C25; Secondary 53A20 }\\
{\bf Keywords:} {Gauss plane, radial density, Four Vertex Theorem}
\vskip 1cm
A manifold with density is a Riemannian manifold $M^n$  with a positive density function $e^{\varphi(x)}$ used to weight volume and hypersurface area. Such manifolds appeared in many ways in mathematics, for example as quotients of Riemannian manifolds or as the Gauss space.  The Gauss space $G^n$ is a Euclidean space with Gaussian probability density $(2\pi)^{-\frac n2}e^{-\frac{r^2}2}$ that is very interesting to probabilists. For more details about manifolds with density and some first results in Morgan's grand goal to ``generalize  all of Riemannian geometry to manifolds with density" we refer the reader to \cite{mo1}, \cite {RCBM}, \cite{Co2}, \cite{ACDLV}.  Following Gromov (\cite[p. 213]{gr1})  the natural generalization of the mean curvature of  hypersufaces on a manifold with density $e^{\varphi}$ is defined as
\begin{equation}      H_{\varphi}=H-\frac 1{n-1}\frac{d\varphi}{d{\bf n}}\end{equation}
and therefore, the curvature of a curve in a plane with density $e^{\varphi}$ is
\begin{equation}      k_{\varphi}=k-\frac{d\varphi}{d{\bf n}}.\end{equation}

We call  $k_{\varphi}$ the curvature with density or $\varphi-$curvature of the curve.

In this note, we study the Four Vertex Theorem in planes with radial density $e^{\varphi(r)}, $ where $r$ is the distance from the origin.  Curves and the function $\varphi$ are assumed to be of the class $C^3$ and $C^2$, respectively.

It is well known that ``\emph{every simple closed curve in the Euclidean plane has at least four vertices}" (the Four Vertex Theorem).  This theorem has a long and interesting history (see \cite{ca}, \cite{turck}).

First we observe that in general, the Four Vertex Theorem does not hold in planes with density.

In  the Gauss plane ${G}^2$  with density
 $e^\varphi=\frac 1{2\pi}e^{-r^2/2}, $  let
$$\begin{aligned}
   \alpha:[0,2\pi]&\longrightarrow\Bbb R^2\\
   t&\longmapsto&(\cos t, \sin t+b)
         \end{aligned}$$
be a parametrization by arc length of the unit circle with center $I(0, b).$
In  the Gauss plane ${G}^2,$ direct computation shows that
\begin{equation}\label{1}
  k_\varphi={x'y''-x''y'}-xy'+x'y. \end{equation}

Applying (\ref{1}), we get the $\varphi$-curvature of the circle
$$k_\varphi=-b\sin t.$$
The equation $k'_\varphi=0$ has exactly two solutions $\frac {\pi}2$ and  $\frac {3\pi}2$ and hence the circle has exactly two vertices.

\begin{figure}
 \begin{center}
  \includegraphics[width=10cm]{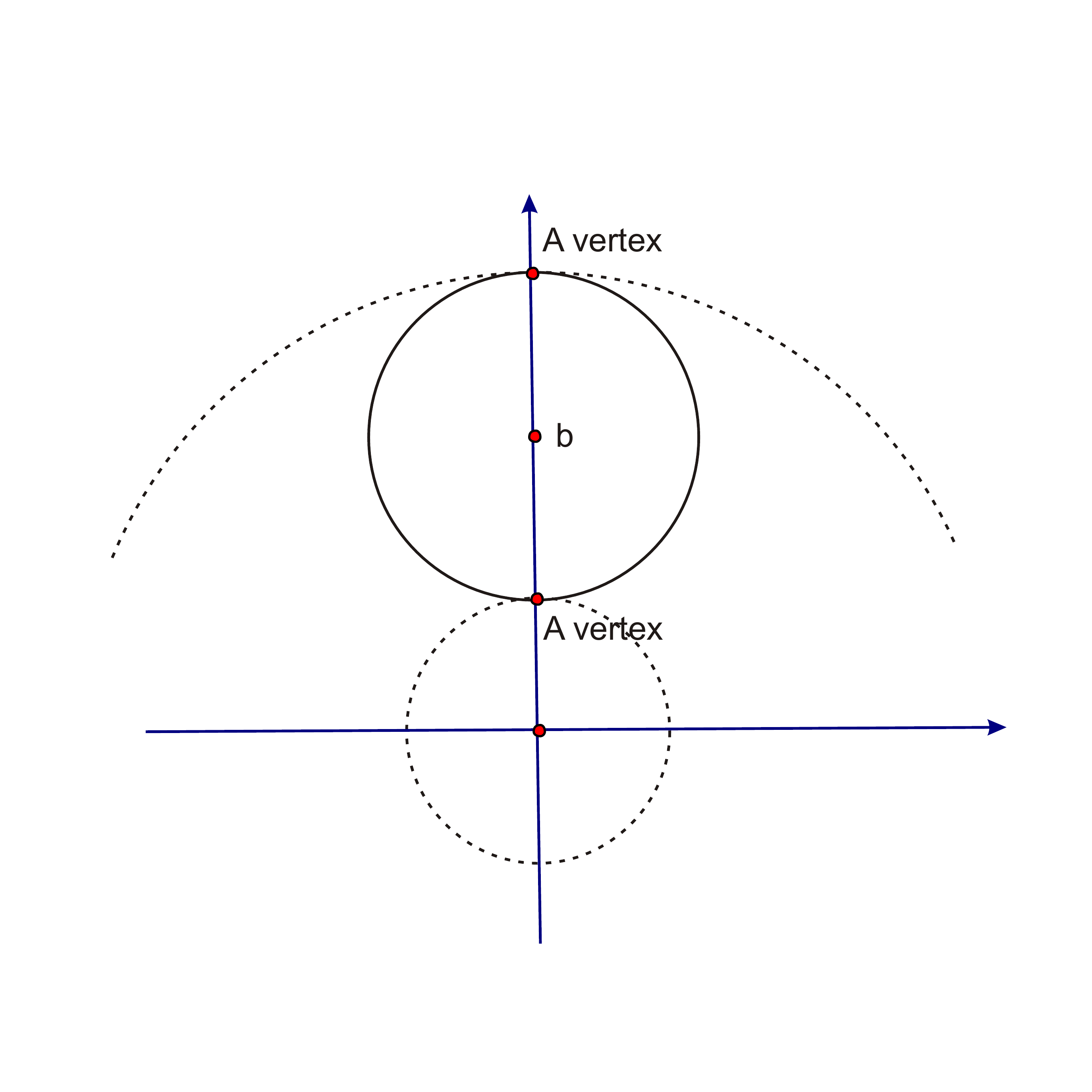}\\
  \end{center}
 \caption{A circle has 2 vertices.}\label{h0}
\end{figure}


The above calculation naturally raises  two questions.
\begin{enumerate}  \item For which classes of simple closed curves the Four Vertex Theorem holds in planes with arbitrary radial densities.
\item Does there exist a plane with radial density other than the Euclidean one, such that the Four Vertex Theorem holds for the class of all closed simple curves?
        \end{enumerate}

We can see that under a rotation about the origin the curvature $k$ of a curve and $\frac {d\varphi}{d{\bf n}}$ are not changed. Thus, the $\varphi$-curvature of a curve in plane with radial density  $e^{\varphi(r)}$ is invariant under a rotation about the origin.
Therefore, as corollaries of this fact, we have the Four Vertex Theorem for the class of simple closed curves that are symmetric about the origin or more generally, for the class of simple closed curves that are invariant under a rotation about the origin.

  \begin{cor} \label{dx} In planes with density
    $e^{\varphi(r)}$, every simple closed curve $\alpha$ symmetric about the origin has at least four vertices.   \end{cor}
\noindent{\bf Proof.}
 Since $k_{\varphi}$ is continuous on a closed curve it reaches a minimum and a maximum on the curve. Suppose that $k_\varphi$  reaches its minimum at $t_1$  and its maximum at $t_2.$ Obviously,  $k_{\varphi}(t_1)=k_{\varphi}(t_1+\pi)$ and  $k_{\varphi}(t_2)=k_{\varphi}(t_2+\pi).$ But we have $\alpha(t_1)\ne \alpha(t_2+\pi)$ and $\alpha(t_2)\ne \alpha(t_1+\pi),$ unless $k_{\varphi}=\text{const.}.$ \hfill $\Box$

By a similar proof, we have a more general result.
\begin{cor}  In planes with density    $e^{\varphi(r)}$, every simple closed curve that is invariant under a rotation about the origin by an angle $\theta< 2\pi$ has at least four vertices.   \end{cor}

For the second question, we have the following as the main result of this note.
\begin{thm} \label{theomain}In a plane with density $e^{\varphi(r)}$, the Four Vertex Theorem holds for the class of all simple closed curves  if and only if $\varphi$ is a constant.\end{thm}

In order to prove Theorem \ref{theomain} we need the following.
\begin{lem} \label{lemmain}
Let $C$ be the circle with center at $(0,b)$ and radius $R$ in the plane $\Bbb R^2$ with density $e^{\varphi(r)}$ and
$    \alpha  : [0,2\pi] \longrightarrow \Bbb R^2;\ \
\alpha(t)=(R\cos t, R\sin t+b)$ be its parametrization. Assume that $C$ does not pass through the origin.
Denote by $k_{\varphi}$ the $\varphi$-curvature of $C.$
 Then we have
 \begin{equation}\label{main}    \frac {dk_{\varphi}}{dt}=\frac{b\cos t}{r^3}\left(\varphi''\frac{r^3+(R^2-b^2)r}2+\varphi' \frac{r^2-(R^2-b^2)}2\right).   \end{equation}
\end{lem}
{\bf Proof.} Direct computation yields
$$k_{\varphi}=\frac 1R+\varphi' \frac{R+b\sin t}{r}.$$
Therefore,
$$  \begin{aligned}
\frac{dk_{\varphi}}{dt}&=\frac{d}{dt}\left(\frac 1R+\varphi'\frac{R+b\sin t}{r}\right) \\
&=\varphi''\frac{dr}{dt}\frac{R+b\sin t}{r}+\varphi'\frac 1{r^2}\left(br\cos t-\frac{dr}{dt}(R+b\sin t)\right)\\
&=\varphi''\frac{bR\cos t}{r}\frac{R+b\sin t}{r}+\varphi'\frac 1{r^2}\left(br\cos t-\frac{bR\cos t}{r}(R+b\sin t)\right)\\
&=\frac{b\cos t}{r^3}\left[\varphi''(R^2r+bRr\sin t)+\varphi'(b^2+bR\sin t)\right].
               \end{aligned}$$
Since on the circle we have $r=\sqrt{R^2+b^2+2Rb\sin t},$ replacing $bR\sin t$ by $\frac{r^2-R^2-b^2}2,$ we obtain  (\ref{main}). \hfill $\Box$

Lemma \ref{lemmain} has some useful corollaries.
\begin{cor}\label{constant} In planes with density $e^{\varphi(r)}$, a circle about the origin has constant $\varphi$-curvature.
\end{cor}
{\bf Proof.} This is immediate since $b=0.$ \hfill $\Box$
\begin{cor} For every positive integer $n\ge 1$, there exists a radial density $e^{\varphi(r)}$ in the plane such that a given circle containing the origin in its interior, but not as the center has exactly $2n$ vertices.
\end{cor}
{\bf Proof.} We can assume that the center of the circle is $(0,b),\ b>0 $ and the radius is $R>0.$ Since the circle contains the origin, we have $R>b$ and therefore $r^2+R^2-b^2>0.$ Let $r_1, r_2,\ldots, r_{n-1}$ be positive real numbers such that $R-b<r_1<r_2<\ldots<r_{n-1}<R+b.$ Set
$$p(r)=(r-r_1)(r-r_2)\ldots(r-r_{n-1}).$$
It is easy to see the ODE
$$\varphi''(r^3+(R^2-b^2)r)+\varphi'(r^2-R^2+b^2)=r^2p(r)$$
has a family of solutions of the following form
$$\varphi=\int\frac{\left(\int p(r)dr+c_1\right)r}{r^2+R^2-b^2}dr+c_2.$$
Therefore, the equation $\frac {dk_{\varphi}}{dt}=0$ has exactly $2n$ solutions and it follows that the circle has exactly $2n$ vertices.
 \hfill $\Box$

\begin{figure}
 \begin{center}  \includegraphics[width=12cm]{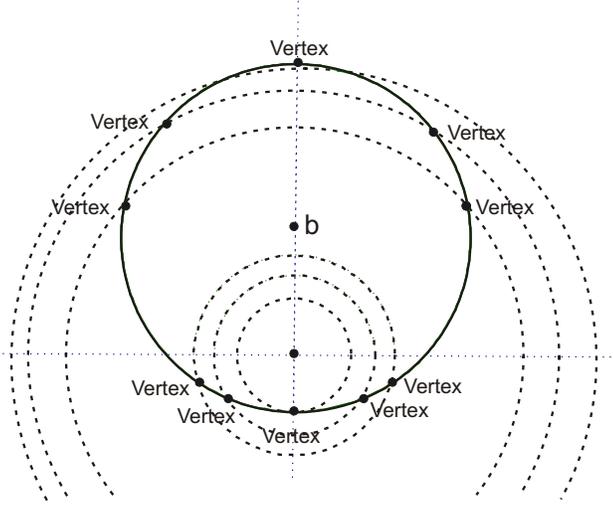}\\\end{center}
  \caption{There exists a radial density such that a given circle containing the origin has $2n$ vertices.}\label{h1}
\end{figure}



\begin{cor} There exists a radial density $e^{\varphi(r)}$ in the plane such that a given circle containing the origin in its interior has constant density curvature.
\end{cor}
{\bf Proof.}  We can assume that the center of the circle is $(0,b),\ b>0 $ (the case of $b=0$ follows from Corollary \ref{constant}) and the radius is $R>0.$ Since the circle contains the origin, we have $R>b$ and therefore $r^2+R^2-b^2>0.$
It is easy to see the ODE
$$\varphi''(r^3+(R^2-b^2)r)+\varphi'(r^2-R^2+b^2)=0$$
has a family of solutions of the following form
$$\varphi=c_1\ln(r^2+R^2-b^2)+c_2.$$
 \hfill $\Box$

\noindent{\bf Proof of Theorem \ref{theomain}.}  We  prove Theorem \ref{theomain} by showing that if $\varphi$ is not a constant, then there exists a circle with center $(0,b)$ and radius $R$ that has exactly two vertices. Now suppose  $\varphi'(b)\ne 0,  b\ne 0.$
\begin{enumerate}
\item If $\varphi'(b) >0,$ then there exists $0<\epsilon<b$ such that  $\varphi'(r) >0, \forall r\in[b-\epsilon, b+\epsilon].$ Let $N=\underset{ r\in[b-\epsilon, b+\epsilon]}\min\varphi'(r) >0.$  Since $\varphi''$ is continuous in $[b-\epsilon, b+\epsilon]$ there exists a positive real number $M>|\varphi''(r)|,\ \forall  r\in [b-\epsilon, b+\epsilon].$

Choose $R>0$ so that $0<R<\min\left\{\epsilon,\frac{b(b-\epsilon)N}{M(b+\epsilon)^2}\right\}.$
We have
$$\begin{aligned}  |\varphi'' Rr(R+b\sin t)|&<MR(b+\epsilon)|R+b\sin t|\\
&\le MR(b+\epsilon)(R+b)\\
&< MR(b+\epsilon)^2\\
&< M\frac{b(b-\epsilon)(b+\epsilon)^2N}{M(b+\epsilon)^2}\\
&=b(b-\epsilon)N\\
&<b(b-R)N\\
&\le \varphi'b(b+R\sin t).\\
                 \end{aligned}$$

Recall that $\varphi'' Rr(R+b\sin t)=\varphi''\frac{r^3+(R^2-b^2)r}2$ and  $\varphi'b(b+R\sin t)=\varphi'\frac{r^2-(R^2-b^2)}2$ (see the proof of Lemma \ref{lemmain}); we have
$$\varphi''\frac{r^3+(R^2-b^2)r}2+\varphi'\frac{r^2-(R^2-b^2)}2>0$$
and hence the circle with center $(0,b)$ and radius $R$ has exactly two vertices.

\item If $\varphi'(b) <0,$ then there exists $0<\epsilon<b$ such that  $\varphi'(r) <0, \forall r\in[b-\epsilon, b+\epsilon].$ Let $N=\underset{ r\in[b-\epsilon, b+\epsilon]}\min |\varphi'(r)| >0.$  We have $\varphi'\le -N, \forall r\in[b-\epsilon, b+\epsilon].$ Let $M, R$ be as in the proof of the first part  and by a similar proof, we obtain
    $$ |\varphi'' Rr(R+b\sin t)|<b(b-R)N\le-b(b+R\sin t)\varphi'.$$
Therefore,
$$\varphi'' Rr(R+b\sin t)+b(b+R\sin t)\varphi'<0,$$
or
$$\varphi''\frac{r^3+(R^2-b^2)r}2+\varphi'\frac{r^2-(R^2-b^2)}2<0.$$
Thus, the circle with center $(0,b)$ and radius $R$ has exactly two vertices.
\end{enumerate}
 \hfill $\Box$
 \vskip.5cm
{\bf Acknkowledgements.} We thank  Professor Frank Morgan for bringing the category of manifolds with density to our attention. We also thank the referee for valuable
comments and suggestions.

\end{document}